\documentclass[12pt]{amsart}
\usepackage{amssymb,amsmath,a4wide}

\newtheorem{theorem}{Theorem}[section]
\newtheorem{lemma}[theorem]{Lemma}

\newtheorem{defn}[theorem]{Definition}
\newtheorem{remark}{\it Remark\/}

\newcommand{\eps}{\varepsilon}
\newcommand{\To}{\longrightarrow}
\newcommand{\pa}{\partial}
\newcommand{\RR}{{\mathbb{R}}}
\renewcommand{\div}{\mathop{{\rm div}}}

\begin{document}

\title[Boundary blow-up]{Boundary blow-up and degenerate equations}%
\author{Satyanad Kichenassamy}%
\address{Laboratoire de Math\'ematiques (UMR 6056), CNRS \&\ Universit\'e
de Reims Champagne-Ardenne, Moulin de la Housse, B.P. 1039, F-51687
Reims Cedex 2\\France}%
\email{satyanad.kichenassamy@univ-reims.fr}%

\thanks{}%

\vskip 1em

\hfill{\bf Appeared in:} \emph{Journal of Functional Analysis},
{\bf 215} : 2 (2004) 271--289.

\vskip 1em
\begin{abstract}
Let $\Omega\subset\RR^2$ be a bounded domain of class
$C^{2+\alpha}$, $0<\alpha<1$. We show that if $u$ is the solution
of $\Delta u = 4\exp(2u)$ which tends to $+\infty$ as
$(x,y)\to\partial\Omega$, then the hyperbolic radius $v=\exp(-u)$
is also of class $C^{2+\alpha}$ up to the boundary. The proof
relies on new Schauder estimates for degenerate elliptic equations
of Fuchsian type.
\end{abstract}
\maketitle

\section{Introduction}

Let $\Omega\subset\RR^2$ be a bounded domain of class
$C^{2+\alpha}$, with $0<\alpha<1$. Consider the Liouville equation
\begin{equation}\label{eq:liouville}
-\Delta u +4e^{2u}=0.
\end{equation}
Let $\mathcal{U}$ be the set of solutions of (\ref{eq:liouville})
which belong to $C^{2+\alpha}(\overline{\Omega})$, and consider
\[ u_\Omega=\sup_{u\in\mathcal{U}}u.
\]
It is known (\emph{see} \cite{B,LN}, the survey \cite{BF}
and its references) that
\begin{enumerate}
\item $u_\Omega$ is finite, positive and smooth in $\Omega$.
\item $u_\Omega$ is the limit of the sequence $(u_n)$ of solutions of
(\ref{eq:liouville}) equal to $n$ on $\pa\Omega$, as $n\to\infty$;
it is called the the \emph{maximal solution} of
(\ref{eq:liouville}) on $\Omega$, and dominates all solutions,
thus providing a universal bound on any classical solution of
(\ref{eq:liouville}), independent of its boundary data.
\item If $\Omega'\subset \Omega$, then $u_{\Omega}\leq u_{\Omega'}$ in
         $\Omega'$.
\item If $\Omega$ is simply connected, one can recover a Riemann
map for $\Omega$ from the \emph{hyperbolic radius}
\[ v_\Omega:=\exp(-u_\Omega).
\]
\item The metric $v_\Omega^{-2}(dx^2+dy^2)$ on $\Omega$ has
constant negative curvature; it generalizes the hyperbolic
(Poincar\'e) metric on the unit disk.
\item Denote by $d(x,y)$ the distance of $(x,y)$ to the boundary. It
is of class $C^{2+\alpha}$ near the boundary. As $d\to 0$,
\[ |u_\Omega/\ln(2d)+1|=O(d)
\]
\cite[th.~4]{BM}, $v_\Omega=2d+o(d)$ \cite[p.~204]{BF}, and
$|\nabla v_\Omega|\to 2$ \cite[th.~3.3]{BE}.
\item If $\Omega$ is (convex and) of class $C^{4+\alpha}$, then
$v_\Omega\in C^{2+\beta}(\overline{\Omega})$ for some $\beta>0$
\cite[th.~2.4]{CF}.
\end{enumerate}
An accurate knowledge of the boundary behavior of $v_\Omega$ has
two applications:
\begin{enumerate}
  \item The actual numerical computation of $u_\Omega$ rests on
  the solution of the Dirichlet problem for (\ref{eq:liouville})
  on the domain $\{d>h\}$, where $h$ is small, taking for
  Dirichlet data the beginning of the expansion of $u_\Omega$
  \cite[section 3.3]{BF}. The better this expansion is known,
  the more accurate the computation.
  \item If $v_\Omega$ admits an expansion to second order, one
  finds formally
  \begin{equation}\label{eq:exp}
      v_\Omega=2d-d^2(\kappa+o(1)),
  \end{equation}
  where $\kappa$ is the curvature of $\pa\Omega$. It follows that
  $u$ is convex near any boundary point at which $\kappa>0$. This
  is a local result which does not require $\Omega$ to be convex
  as a whole. This computation is justified by remarks
  \ref{rem:kappa} and \ref{rem:kappa2}, in section \ref{sec:e4}
  of the present paper.
\end{enumerate}
There is a extensive literature on the issue of boundary blow-up;
see \cite{BE,BF,BM,K,KN,LM,LN,MV,O,P} and their references for
further details.

We prove in this paper:
\begin{theorem}\label{th:e1}
If $\Omega$ is of class $C^{2+\alpha}$, then $v_\Omega$ is of
class $C^{2+\alpha}$ near and up to the boundary.
\end{theorem}
\begin{remark} Theorem \ref{th:e1} was conjectured
in \cite[p.~204]{BF}. The result is optimal since $\kappa$ is
precisely of class $C^\alpha$, and not better in general.
\end{remark}
\begin{remark}
Since $u_\Omega$ is smooth inside $\Omega$, we need only
investigate its boundary behavior. Interior bounds on $u_\Omega$
may be obtained by comparison with the exact solutions on balls
containing, or contained in $\Omega$.
\end{remark}
\begin{remark}
From theorem \ref{th:e1}, it follows that $v_\Omega$ solves
\[v_\Omega\Delta v_\Omega=|\nabla v_\Omega|^2-4
\]
up to $\pa\Omega$ in the classical sense.
\end{remark}

An outline of the proof is presented in the following section.
It uses auxiliary results proved in the other sections of the paper.


\section{Outline of proof and organization of the paper}

The procedure consists in reducing the problem to a regularity
problem for a degenerate equation of Fuchsian type, and to prove
estimates which play the role of the boundary Schauder estimates
for the Laplacian. The Fuchsian form shall also make it easy to
find new sub- and super-solutions.

The reduction of a nonlinear PDE to Fuchsian form (see \cite{SK}
and its references) has been useful for constructing explosive
solutions for problems of hyperbolic type; we adapt it to elliptic
problems: for the problem at hand, let us define the
``renormalized unknown'' $w$ by
\[
v_\Omega=2d+d^2w(x,y).
\]
This new unknown solves, near the boundary, the nonlinear Fuchsian
equation
\begin{equation}\label{eq:fuchs}
Lw+2\Delta d = \frac{d^2}{2+dw}\left[2w\nabla w\cdot\nabla d +
d|\nabla w|^2 \right]-2dw\Delta d,
\end{equation}
where
\begin{equation}\label{eq:L}
L=\div(d^2 \nabla)-2= d^2\Delta +2d\nabla d\cdot\nabla-2.
\end{equation}
Recall that an elliptic equation is Fuchsian if (i) its
characteristic form, divided by $d^2$, is uniformly positive
definite; (ii) the first-order terms are $O(d)$ and (iii) the
terms of order zero are bounded near the boundary. There is a
sizable literature on weighted Schauder estimates for elliptic
problems, see \cite{GT,GS} for instance.

Equation (\ref{eq:fuchs}) needs only to be studied in the
neighborhood of the boundary. Let us therefore introduce a
$C^{2+\alpha}$ thin domain $\Omega'\subset \Omega$, on which $d$
is of class $C^{2+\alpha}$ and does not exceed $\delta\leq 1/2$,
such that $\pa\Omega'$ consists of two portions of class
$C^{2+\alpha}$, of which one is $\pa\Omega$ and the other will be
called $\Gamma$.

Equation (\ref{eq:fuchs}) may be rewritten as a linear equation
with $w$-dependent coefficients: for any $f$, we define
\begin{equation}\label{eq:M}
M_w(f) = \frac{d^2}{2+dw}\left[2f\nabla w\cdot\nabla d + d\nabla
w\cdot\nabla f \right]-2df\Delta d.
\end{equation}
We therefore have
\[ (L-M_w)w+2\Delta d=0.
\]

A comparison argument, similar to the one in \cite{BE} for
instance, yields
\begin{theorem}\label{th:e2}
$w$ et $d^2\nabla w$ are bounded near $\pa\Omega$.
\end{theorem}
This theorem is proved in section \ref{sec:e2}. It provides just
enough regularity on the coefficients of $L-M_w$ to put it within
the scope of the analogue, for the operators at hand, of the
$C^{1+\alpha}$ estimate for elliptic operators (theorem
\ref{th:e7}, proved in section \ref{sec:e7}). We apply this result
in section \ref{sec:e3}, and obtain
\begin{theorem}\label{th:e3}
If $\delta$ is small, $dw$ and $d^2\nabla w$ belong to
$C^\alpha(\overline{\Omega'})$, and $d\nabla w$ is bounded near
$\pa\Omega$.
\end{theorem}
Next, one subtracts from $w$ a function $w_0$ such that $w-w_0$ is
sufficiently flat, and which has the regularity we expect $w$ to
have. The function $w_0$ is constructed in section \ref{sec:e4};
the result is:
\begin{theorem} \label{th:e4}
If $\delta$ is small, there is a function $w_0$ such that $w_0$,
$d\nabla w_0$, and $d^2\nabla^2 w_0$ belong to
$C^\alpha(\overline{\Omega'})$, and
\[ Lw_0+2\Delta d=0\]
near $\pa\Omega$.
\end{theorem}
It follows that $d^2w_0$ is of class $C^{2+\alpha}$ near the
boundary. Letting $\tilde{w}=w-w_0$, we construct sub- and
super-solutions which show (section \ref{sec:e5}) that
\begin{theorem}\label{th:e5} There is a constant $\gamma$ such
that
\[ |\tilde{w}|\leq \gamma d\ln (1/d)\]
near $\pa\Omega$.
\end{theorem}
Thanks to a sharpened $C^{1+\alpha}$ estimate (theorem
\ref{th:e8}, proved in section \ref{sec:e8}), one also proves that
$L\tilde{w}$ is equal to the product of $d$ by a function of class
$C^\alpha(\overline{\Omega'})$. Using then a scaled $C^{2+\alpha}$
estimate (theorem \ref{th:e9}, proved in section \ref{sec:e9}), it
follows (section \ref{sec:e6}) that
\begin{theorem}\label{th:e6}
If $\delta$ is small, $d^2\tilde{w}\in
C^{2+\alpha}(\overline{\Omega})$.
\end{theorem}
Theorem \ref{th:e1} follows.

Section \ref{sec:prel} collects basic notation and computations
which will be used in the paper. Section \ref{sec:e2} gives the
first comparison argument, proving theorem \ref{th:e2}, and
section \ref{sec:e5} the sub- and super-solution argument showing
that $\tilde{w}$ is flat near the boundary. Section \ref{sec:e4}
gives the construction of $w_0$. All general-purpose Schauder-type
estimates are collected in section \ref{sec:fuchsian}.


\section{Preliminary computations}
\label{sec:prel}

We collect simple formulae which will be useful in the sequel, and
which follow by direct computation. Fix a point $P$ on
$\pa\Omega$, which we take as origin of coordinates in $\RR^2$;
define the change of variables $(x,y)\mapsto (T,Y)$, where
\[ T=d(x,y) \text{ and } Y=y. \]
It is well-defined near the boundary, and of class $C^{2+\alpha}$.
We may also assume, by performing a rigid motion, that $\pa d/\pa
x=1$ and $\pa d/\pa y=0$ at $P$; the $y$-axis is then tangent to
the boundary at $P$. The Jacobian of the change of variables is
$d_x$, which equals $1$ at $P$; the change of variables is
therefore invertible, and of class $C^{2+\alpha}$ together with
its inverse, if $(x,y)$ is small.

If $\kappa$ denotes the curvature of the boundary, and subscripts
denote derivatives,
\begin{align}
|\nabla d|&=1, \quad \Delta d=-\frac\kappa{1-T\kappa};  \\ \pa_x
&= d_x\pa_T, \quad \pa_y = d_y\pa_T+\pa_Y, \quad d_y=d_Y;\\ \Delta
w&=w_{TT}+w_{YY}+2d_yw_{TY}+w_T\Delta d.
\end{align}
Let
\[ D=T\pa_T, \quad \Delta'=\pa_{TT}+\pa_{YY}. \]
We find
\begin{equation}
e^{-u}[-\Delta u+4e^{2u}] = Lw+2\Delta d-M_w(w), \label{eq:ssol}
\end{equation}
with $L$ given by equation (\ref{eq:L}) and $M_w$ by (\ref{eq:M}),
and
\begin{align}
\nabla d\cdot\nabla w
  &= d_xw_x+d_yw_y= d_x^2w_T+d_y(d_yw_T+w_Y)\nonumber\\
  &= w_T+d_yw_Y;\nonumber\\
Lw&=d^2\Delta w+2d\nabla d\cdot\nabla w-2w \nonumber\\
  &= T^2\Delta w+2T(w_T+d_yw_Y)-2w\nonumber\\
  &= T^2(\Delta'w+2d_yw_{TY}+w_T\Delta d)+(2D-2)w+2Td_yw_Y\nonumber\\
  &= T^2w_{TT}+2(D-1)w+T(\Delta d)w_T+2Td_y\pa_Y(D+1)w;\nonumber\\
L &=L_0+L_1,\nonumber\\ L_0&=(D+2)(D-1)+T^2\pa^2_Y,\\
L_1&=2Td_y(D+1)\pa_Y+T(\Delta d)D.
\end{align}

We also need the spaces $C^{k+\alpha}_\sharp(U)$, for $k=1$ or 2,
and any $U\subset \Omega$:
\begin{defn}
We say that $u\in C^{k+\alpha}_\sharp(U)$ if $T^j u\in
C^{j+\alpha}(U)$ for $0\leq j\leq k$. Its norm is the sum of the
$\|T^j u\|_{C^{j+\alpha}}$.
\end{defn}
It is equivalent to require that $T^j \nabla^j u\in C^{\alpha}(U)$
for $0\leq j\leq k$.

There are two auxiliary domains which will be used for
localization.

The first is the domain $\Omega'\subset\Omega$ already mentioned,
which is such that $\pa\Omega'=\pa\Omega\cup\Gamma$.

The second is defined in the $(T,Y)$ coordinates, by
\[
\Omega''=\{(T,Y) : 0<T<\theta, |Y|<\theta\},
\]
where $\theta$ will be chosen small in \ref{ssec:sol-w0}. Note
that since $d_y(P)=0$, it is $O(\theta)$ over $\Omega''$, and
therefore
\begin{equation}\label{eq:pert}
\|L_1w\|_{C^{\alpha}(\overline{\Omega''})}\leq
c\theta\|w\|_{C^{2+\alpha}_\sharp(\Omega'')},
\end{equation}
where $c$ is independent of $\theta$.


\section{Proof of theorem \ref{th:e2}}
\label{sec:e2}

By comparison with the maximal solution on balls entirely
contained in $\Omega$, we obtain interior bounds. It suffices to
find bounds near the boundary. We write $u$ instead of $u_\Omega$,
for short.

$\pa\Omega$ satisfies a uniform interior and exterior sphere
condition at every point. Furthermore, there is an $r_0>0$ such
that any point $P$ such that $d(P)< r_0$ admits a unique nearest
point $Q$ on the boundary. Making $r_0$ smaller if necessary, we
may assume that there are two points $A$ and $A'$ such that the
balls $B_{r_0}(A)$ and $B_{r_0}(A')$ are tangent to $\pa\Omega$ at
$Q$ and furthermore
\[ \Omega_i \subset \Omega\subset\Omega_e,  \]
where $\Omega_i=B_{r_0}(A)$ and $\Omega_e=B_{1/r_0}(A')\setminus
B_{r_0}(A')$. The line segment $AQ$ is a radius of $B_{r_0}(A)$.

Let $u_e$ and $u_i$ be the maximal solutions of
(\ref{eq:liouville}) on $\Omega_e$ and $\Omega_i$ respectively.
They are known explicitly: they are radial, and satisfy the
conclusion of theorem \ref{th:e1}. Therefore,
$\exp(-u_e)=d(2+dw_e)$ and $\exp(-u_i)=d(2+dw_i)$, where $w_e$ and
$w_i$ are bounded over $AQ$, by a quantity which depends only on
$r_0$, and not on $A$.
\begin{remark}
In fact, for any point $P$, if $r=AP$, $r'=A'P$, we have
$v_i(P)=r_0-r^2/r_0$, and $v_e(P)=4\pi^{-1}\ln
r_0\cos(\frac\pi{2\ln r_0}\ln r')r'$ \cite[p.~201]{BF}.
\end{remark}
Since solutions to (\ref{eq:liouville}) decrease as $\Omega$
increases, we have
\[ u_e\leq u\leq u_i\text{ over } B_{r_0}(A).\]

Therefore, $w$ is bounded over the segment $AQ$.

In particular, $|w|$ is bounded over $\{P : d(P)<r_0\}$ by some
number $M$, since $P$ lies on the corresponding segment $AQ$.
Therefore,
\[ 2d-Md^2\leq \exp (-u)\leq 2d+Md^2.
\]

We now use scaling and regularity estimates (as in
\cite[th.~3.3]{BE},\cite[lemma 2.2, p.~289]{SK-th}) to derive
gradient bounds from pointwise bounds. Consider $P$ such that
$d(P)=2\sigma$ with $3\sigma<r_0$. For $(x,y)$ in the unit disk,
let
\[P_\sigma=P+(\sigma x,\sigma y)\]
and
\[ u_\sigma(x,y):=u(P_\sigma)+\ln \sigma.
\]
One verifies that $u_\sigma$ solves (\ref{eq:liouville}).

Since $\sigma<d(P_\sigma)<3\sigma$, we have, for $r_0$ so small
that $2d\pm Md^2$ is an increasing function of $d$ for $d<r_0$,
\[ 2\sigma-M\sigma^2 < \exp(- u(P_\sigma))=\sigma\exp(-u_\sigma(x,y)) <
6\sigma+9M\sigma^2,
\]
hence
\[ 2-M\sigma <\exp(-u_\sigma(x,y))< 6+9M\sigma.
\]
It follows that $\exp(-u_\sigma)$ is bounded and bounded away from
zero on the unit ball if $\sigma$ is small. It follows that
$u_\sigma$ itself is bounded. By interior regularity, it is
bounded in $C^1$ on the ball of radius one-half. Applying this
result at the origin, we find, recalling that $\sigma=\frac12
d(P)$,
\[ u(P)+\ln d(P) \text{ and } d\nabla u(P) \text{ are bounded near }\pa\Omega.
\]
Since $u=-\ln(2d+d^2w)=-\ln d-\ln(2+dw)$,
\[ d\nabla u = -\nabla d -(2+dw)^{-1}d[w\nabla d+d\nabla w],\]
and since $|\nabla d|=1$, and we already know that $w$ is bounded,
we find that
\[ w(P) \text{ and } d^2\nabla w(P) \text{ are bounded near }\pa\Omega,
\]
QED.


\section{Two types of Fuchsian operators}
\label{sec:fuchsian}

\subsection{Scaled Schauder estimates}

Theorems \ref{th:e3} and \ref{th:e6} follow from general
Schauder estimates for linear Fuchsian operators, applied to
$L-M_w$. We need to distinguish two types of operators, according
to the regularity of their coefficients.

An operator $A$ is said to be \emph{of type (I)} (on a given
domain) if it can be written
\[ A=\pa_i(d^2a^{ij}\pa_{j})+db^i\pa_i+c, \]
with $(a^{ij})$ uniformly elliptic and of class $C^\alpha$, and
$b^i$, $c$ bounded.
\begin{remark}
One can also allow terms of the type $\pa_i(b^{\prime i}u)$ in
$Au$, if $b^{\prime i}$ is of class $C^\alpha$, but this
refinement will not be needed here.
\end{remark}

An operator is said to be \emph{of type (II)} if it can be written
\[ A= d^2a^{ij}\pa_{ij}+db^i\pa_i+c, \]
with $(a^{ij})$ uniformly elliptic and $a^{ij}$, $b^i$, $c$ of
class $C^\alpha$.
\begin{remark}
One checks directly that types (I) and (II) are invariant under
changes of coordinates of class $C^{2+\alpha}$. In particular, to
check that an operator is of type (I) or (II), we may work
indifferently in coordinates $(x,y)$ or $(T,Y)$ defined in section
\ref{sec:prel}. All proofs will be performed in the $(T,Y)$
coordinates; an operator is of type (II) precisely if it has the
above form with $d$ replaced by $T$, and the coefficients
$a^{ij}$, $b^i$, $c$ are of class $C^\alpha$ as functions of $T$
and $Y$; a similar statement holds for type (I).
\end{remark}

The basic results are
\begin{theorem}\label{th:e7}
If $Ag=f$, where $f$ et $g$ are bounded and $A$ is of type (I) on
$\Omega'$, then $d\nabla g$ is bounded, and $dg$ and $d^2\nabla g$
belong to $C^\alpha(\Omega'\cup\pa\Omega)$.
\end{theorem}
\begin{theorem}\label{th:e8}
If $Ag=df$, where $f$ and $g$ are bounded, $g=O(d^\alpha)$, and
$A$ is of type (I) on $\Omega'$, then $g\in
C^{\alpha}(\Omega'\cup\pa\Omega)$ and $dg\in
C^{1+\alpha}(\Omega'\cup\pa\Omega)$
\end{theorem}
\begin{theorem}\label{th:e9}
If $Ag=df$, where $f\in C^\alpha(\Omega'\cup\pa\Omega)$,
$g=O(d^\alpha)$, and $A$ is of type (II) on $\Omega'$, then $d^2
g$ belongs to $C^{2+\alpha}(\Omega'\cup\pa\Omega)$.
\end{theorem}

Let $\rho>0$ and $t\leq 1/2$. Throughout the proofs, we shall use
the sets
\begin{align*}
Q   &=\{ (T,Y) : 0\leq T\leq 2 \text{ and } |y|\leq 3\rho\},\\
Q_1 &=\{ (T,Y) : \frac14\leq T\leq2 \text{ and } |y|\leq 2\rho\},\\
Q_2 &=\{ (T,Y) : \frac12\leq T\leq 1 \text{ and } |y|\leq\rho/2\},\\
Q_3 &=\{ (T,Y) : 0\leq T\leq \frac12 \text{ and } |y|\leq\rho/2\}.
\end{align*}
We may assume, by scaling coordinates, that
$Q\subset\Omega'$. It suffices to prove the announced regularity
on $Q_3$.


\subsection{Proof of theorem \ref{th:e7}}
\label{sec:e7}

Let $Af=g$, with $A$, $f$, $g$ satisfying the assumptions of the
theorem over $Q$, and let $y_0$ be such that $|y_0|\leq \rho$.

For $0<\eps\leq 1$, and $(T,Y)\in Q_1$, let
\[ f_\eps(T,Y)=f(\eps T,y_0+\eps Y),
\]
and similarly for $g$ and other functions. We have
$f_\eps=(Ag)_\eps=A_\eps f_\eps$, where
\[
A_\eps=\pa_i(T^2a^{ij}_\eps\pa_{j})+Tb^i_\eps\pa_i+c_\eps
\]
is also of type (I), with coefficient norms independent of
$\eps$ and $y_0$, and is uniformly elliptic in $Q_1$.

Interior estimates give
\begin{equation}\label{eq:c-1-alpha}
\|g_\eps\|_{C^{1+\alpha}(Q_2)}\leq
M_1:=C_1(\|f_\eps\|_{L^{\infty}(Q_1)}+
\|g_\eps\|_{L^{\infty}(Q_1)}).
\end{equation}
The assumptions of the theorem imply that $M_1$ is independent of
$\eps$ and $y_0$.

We therefore find,
\begin{align}
|\eps\nabla g(\eps T, y_0+\eps Y)| &\leq M_1,\\
\eps|\nabla g(\eps T, y_0+\eps Y)-\nabla g(\eps T', y_0)|
&\leq M_1
   (|T-T'|+|Y|)^\alpha\,
\end{align}
if $\frac12\leq T, T'\leq 1$ and $|Y|\leq\rho/2$.
It follows in particular, taking
$Y=0$, $\eps=t\leq 1$, $T=1$, and recalling that $|y_0|\leq\rho$, that
\begin{equation}
  |t\nabla g(t,y)|\leq M_1\text{ if } |y|\leq \rho, t\leq 1.
\end{equation}
This proves the first statement in the theorem.

Taking
$\eps=2t\leq 1$, $T=1/2$, and letting $y=y_0+\eps Y$, $t'=\eps T'$,
\[
  2t|\nabla g(t,y)-\nabla g(t',y_0)|\leq
      M_1(|t-t'|+|y-y_0|)^\alpha (2t)^{-\alpha}
\]
for $|y-y_0|\leq \rho t$ and $t\leq t'\leq 2t\leq 1$.

Let us prove that
\begin{equation}\label{eq:est}
 |t^2\nabla g(t,y)-t^{\prime 2}\nabla g(t',y_0)|\leq
      M_2(|t-t'|+|y-y_0|)^\alpha
\end{equation}
for $|y|, |y_0|\leq \rho$, and $0\leq t\leq t'\leq \frac12$, which
will prove
\[  t^2\nabla g\in C^\alpha(Q_3).
\]
It suffices to prove this estimate in the two cases: (i) $t=t'$
and (ii) $y=y_0$; the result then follows from the triangle
inequality. We distinguish three cases.
\begin{enumerate}
  \item If $t=t'$, we need only consider the case $|y-y_0|\geq\rho
  t$. We then find
  \[ t^2|\nabla g(t,y)-\nabla g(t,y_0)|\leq
      2M_1t\leq 2M_1|y-y_0|/\rho.
  \]
  \item If $y=y_0$ and $t\leq t'\leq 2t\leq 1$, we have $t+t'\leq
  2t'$, hence
  \begin{align*}
      |t^2\nabla g(t,y_0)-t^{\prime 2}\nabla g(t',y_0)|
         &\leq t^2|\nabla g(t,y_0)-\nabla g(t',y_0)|
            +|t-t'|(t+t')|\nabla g(t',y_0)|    \\
         &\leq M_12^{-1-\alpha}t^{1-\alpha}|t-t'|^\alpha
                +2M_1|t-t'|                    \\
         &\leq M_2|t-t'|^\alpha.
  \end{align*}
  \item If $y=y_0$, and $2t\leq t'\leq 1/2$, we have $t+t'\leq
  3(t'-t)$, and
  \begin{align*}
      |t^2\nabla g(t,y_0)-t^{\prime 2}\nabla g(t',y_0)|
         &\leq  M_1(t+t') \\
         &\leq 3M_1 |t-t'|.
  \end{align*}
\end{enumerate}
This proves estimate (\ref{eq:est}).

On the other hand, since $g$ and $T\nabla g$ are bounded over $Q_3$,
\[  Tg\in \mathrm{Lip}(Q_3)\subset C^\alpha(Q_3).
\]
This completes the proof.


\subsection{Proof of theorem \ref{th:e8}}
\label{sec:e8}

The argument is similar, except that $M_1$ is now replaced by
$M_3\eps^\alpha$, with $M_3$ independent of $\eps$ and $y_0$.
It follows that
\begin{equation}
  |t\nabla g(t,y)|\leq M_3t^\alpha\text{ if } |y|\leq \rho, t\leq 1.
\end{equation}
Taking $\eps=2t\leq 1$, $T=1/2$, and letting $y=y_0+\eps Y$,
$t'=\eps T'$, and noting that
$\eps^\alpha(|T-T'|+|Y|)^\alpha=(|t-t'|+|y-y_0|)^\alpha$, we find
\[
  2t|\nabla g(t,y)-\nabla g(t',y_0)|\leq
      M_3 (|t-t'|+|y-y_0|)^\alpha
\]
for $|y-y_0|\leq \rho t$ and $t\leq t'\leq 2t\leq 1$.
Let us prove that
\begin{equation}\label{eq:est2}
 |t\nabla g(t,y)-t'\nabla g(t',y_0)|\leq
      M_4(|t-t'|+|y-y_0|)^\alpha
\end{equation}
for $|y|, |y_0|\leq \rho$, and $0\leq t\leq t'\leq \frac12$, which
will prove
\[  T\nabla g\in C^\alpha(Q_3).
\]
We again distinguish three cases.
\begin{enumerate}
  \item If $t=t'$, $|y-y_0|\geq\rho t$, we find
  \[ t|\nabla g(t,y)-\nabla g(t,y_0)|\leq
      2M_3t^\alpha\leq 2M_3(|y-y_0|/\rho)^\alpha.
  \]
  \item If $y=y_0$ and $t\leq t'\leq 2t\leq 1$, we have $|t-t'|\leq
  t\leq t'$, hence
  \begin{align*}
      |t\nabla g(t,y_0)-t'\nabla g(t',y_0)|
         &\leq \frac12 M_3|t-t'|^\alpha
                +|t-t'||\nabla g(t',y_0)|             \\
         &\leq M_3|t-t'|^\alpha(\frac12+t'^{1-\alpha}t'^{\alpha-1})
               \leq 2M_3|t-t'|^\alpha.
  \end{align*}
  \item If $y=y_0$, and $2t\leq t'\leq 1/2$, we have $t\leq t'\leq
  3(t'-t)$, and
  \begin{align*}
      |t\nabla g(t,y_0)-t'\nabla g(t',y_0)|
         &\leq  M_3(t^\alpha+t^{\prime \alpha}) \\
         &\leq 2M_3 (3|t-t'|)^\alpha.
  \end{align*}
\end{enumerate}
Estimate (\ref{eq:est2}) therefore holds.

The same type of argument shows that
\[  g\in C^\alpha(Q_3). \]
In fact, we have, with again $\eps=2t$,
$\|g_\eps\|_{C^\alpha(Q_2)}\leq M_5\eps^\alpha$, where $M_5$
depends on the r.h.s.\ and the uniform bound assumed on $f$. This
implies
\[ |g(t,y)-g(t',y_0)|\leq M_5 (|t-t'|+|y-y_0|)^\alpha,
\]
if  $t\leq t' \leq 2t\leq 1\text{ and } |y-y_0|\leq\rho t$. The
assumptions of the theorem yield in particular
\[ |g(t,y)|\leq M_5 t^\alpha,
\]
for $t\leq 1/2$ and $|y|\leq \rho$.

If $\rho t\leq |y-y_0|\leq \rho$, and $t\leq 1/2$, we have
\[ |g(t,y)-g(t,y_0)|\leq 2M_5t^\alpha
     \leq 2M_5\left(\frac{|y-y_0|}{\rho}\right)^\alpha.
\]

If $2t\leq t'\leq 1/2$ and $y=y_0$,
\[ |g(t,y_0)-g(t',y_0)|
  \leq M_5(t^\alpha+t'^\alpha)
  \leq 2 M_5(3|t-t'|)^\alpha.
\]

If $t\leq t'\leq 2t\leq 1/2$, we already have
\[ |g(t,y_0)-g(t',y_0)| \leq M_5 |t-t'|^\alpha.
\]
The H\"older continuity of $g$ follows.

Combining these pieces of information, we conclude that
\[  g\in C^{1+\alpha}_\sharp(Q_3), \]
QED.

\subsection{Proof of theorem \ref{th:e9}}
\label{sec:e9}

We must now use interior $C^{2+\alpha}$ estimates, rather than
$C^{1+\alpha}$ estimates. We therefore have, instead of equation
(\ref{eq:c-1-alpha}),
\begin{equation}\label{eq:c-2-alpha}
\|g_\eps\|_{C^{2+\alpha}(Q_2)}\leq
C_2(\|g_\eps\|_{L^{\infty}(Q_1)}+ \|f_\eps\|_{C^{\alpha}(Q_1)}).
\end{equation}
The assumptions guarantee that this quantity is $O(\eps^\alpha)$.
The previous argument ensures that $g$ and $d\nabla g$ belong to
$C^\alpha(Q_3)$; furthermore, we also have
\[
|t^2\nabla^2g| \leq M_6 t^\alpha,\quad \text{ for } |y|\leq \rho,
t\leq 1
\]
and
\[ t^{2}|\nabla^2g(t,y)-\nabla^2g(t',y_0)|\leq
  M_6(|t-t'|+|y|)^\alpha,
\]
for
\[ t\leq t'\leq 2t\leq 1\text{ and } |y|\leq \rho t.
\]
\begin{enumerate}
\item If $\rho t\leq |y-y_0|\leq \rho$, and $t\leq 1/2$, we have
  \[ t^2|\nabla^2 g(t,y)-\nabla^2 g(t,y_0)|\leq 2M_6t^\alpha
     \leq 2M_6\left(\frac{|y-y_0|}{\rho}\right)^\alpha.
  \]
\item If $2t\leq t'\leq 1$ and $y=y_0$,
  \[ |t^2\nabla^2 g(t,y_0)-t^{\prime 2}\nabla^2 g(t',y_0)|\leq
   M_6(t^\alpha+t'^\alpha)
  \leq 2M_6(3|t-t'|)^\alpha.
  \]
\item If $t\leq t'\leq 2t\leq 1$, we have
  \begin{align*}
    |t^2\nabla^2 g(t,y_0)-t^{\prime 2}\nabla^2 g(t',y_0)|
   &\leq M_6|t-t'|^\alpha+
    |t-t'|(t+t')|\nabla^2 g(t',y_0)| \\
   &\leq M_6|t-t'|^\alpha+|t-t'|^\alpha
      t^{1-\alpha}(2t')M_6t^{\prime\alpha-2} \\
   &\leq 3M_6|t-t'|^\alpha.
  \end{align*}
\end{enumerate}
It follows that
\[  t^2\nabla^2 g\in C^\alpha(Q_3). \]
By inspection, the second derivatives of $t^2g$ are all of class
$C^\alpha$, taking into account the fact that $g$ and $t\nabla g$
are. We conclude that
\[  T^2g\in C^{2+\alpha}(Q_3), \]
QED.

\subsection{Proof of theorem \ref{th:e3}}
\label{sec:e3}

Since $d$ is $C^{2+\alpha}$, and theorem \ref{th:e2} gives us that
$w$ and $d^2\nabla w$ are bounded, we have near $\pa\Omega$
\begin{enumerate}
  \item operator $L-M_w$ is of type (I);
  \item $(L-M_w)w$ and $w$ are bounded;
\end{enumerate}
theorem \ref{th:e7} therefore applies. The desired conclusion follows.


\subsection{Proof of theorem \ref{th:e6}}
\label{sec:e6}

It suffices to show that $d^2\tilde w$ is of class $C^{2+\alpha}$
near (and up to) $\pa\Omega$.

Equation (\ref{eq:fuchs}) now takes the form
\[ L\tilde{w} =M_w(w),\]
where we know from theorem \ref{th:e5} that
$\tilde{w}=O(d\ln(1/d))$ and from theorem \ref{th:e3}
that $d\nabla w$ is bounded.

Using the expression of $M_w(w)$, we find that
\[ L\tilde{w} \in dL^\infty. \]

Since $L$ is of type (I), theorem \ref{th:e8} now tells us that
$\tilde{w}$ and $d\nabla\tilde{w}$ are of class $C^\alpha$.

Thanks to the regularity of $w_0$, we infer that $w$ and $d\nabla
w$ are $C^\alpha$. We therefore find that in fact,
\[ L\tilde{w} \in dC^\alpha. \]

Since $L$ is also of type (II), theorem \ref{th:e9} now enables us
to conclude that $d^2\tilde{w}\in C^{2+\alpha}$, QED.


\section{Construction of $w_0$ and proof of theorem \ref{th:e4} }
\label{sec:e4}

We localize the problem, and work on the set
$\Omega''=(0,\theta)\times \{|Y|<\theta\}$ associated to a point
$P$ on the boundary, as described in section \ref{sec:prel}.
Recall that, performing a rigid motion if necessary, we may assume
that $\nabla d= (1,0)$ at $P$. One then performs the change of
coordinates $(x,y)\mapsto(T,Y)$, where $T=d(x,y)$ and $Y=y$.

Recall also, from section \ref{sec:prel}, that in coordinates
$(T,Y)$, $L$ takes the form $L=L_0+L_1$, where
\[ L_0 = (D+2)(D-1)+T^2\pa^2_Y.  \]
Furthermore, $\|L_1w\|_{C^{\alpha}(\overline{\Omega''})}\leq
c(\theta)\|w\|_{C^{2+\alpha}_\sharp(\Omega'')}$, where $c(\theta)$
is small if $\theta$ is small. Throughout, we will be only
interested in regularity near $T=Y=0$.

We shall prove that equation
\[
Lw_0=k(T,Y)
\]
admits, for $k\in C^\alpha(\overline{\Omega''})$, such that
$k(T,-\theta)=k(T,\theta)$, a solution in
$C^{2+\alpha}_\sharp(\Omega'')$, which is periodic of period
$2\theta$ with respect to $Y$. Using a partition of unity, it
follows that $Lw_0=k$ admits, near the boundary of $\Omega$, a
solution having the regularity properties required in theorem
\ref{th:e4}.

\subsection{Solution of $L_0w_1=k(T,Y)$.}

Let $F_1 : C^\alpha(\overline{\Omega''})\To C^\alpha(T\geq 0)$
denote a bounded extension operator, such that for any function
$k$, $F_1[k]$ (i) vanishes for $T\geq 2$, (ii) is
$2\theta$-periodic in $Y$ and (iii) coincides with $k$ in
$\Omega''$.

Let
\begin{align*}
F_2 : C^\alpha(\overline{\Omega''})&\To C^\alpha(T\geq 0)\\
k&\mapsto\tilde{k},
\end{align*}
where
\[ \tilde{k}=\int_1^\infty
F_1[k](T\sigma,Y)\frac{d\sigma}{\sigma^2}.
\]
Note that $(D-1)\tilde{k}=-k$. Since $\int_1^\infty
\sigma^{\alpha-2}d\sigma<\infty$, one checks that $\tilde{k}$ is
indeed in $C^\alpha(T\geq 0)$. We also have, for $T=0$,
$\tilde{k}(0,Y)=k(0,Y)$. Since $D\tilde{k}=\tilde{k}-k$, we find
that $D\tilde{k}$ also is of class $C^\alpha$.

Next, find $h(x,y)$ by solving
\[\Delta' h+\tilde{k}=0,
\]
with $h=0$ for $T=0$, $h_T=0$ for $T=\theta$, and periodic
boundary conditions in $Y$: $h(T,Y+2\theta)=h(T,Y)$; $h$ is
therefore in $C^{2+\alpha}(\overline{\Omega''})$, by the usual
Schauder theory. In particular, $Dh=Th_T=0$ for $T=0$ and
$T=\theta$.

Finally, let
\[w_1=T^{-2}[(D-1)h].
\]
Note that $w_1=T^{-1}D(h/T)$.
\begin{remark}\label{rem:kappa}
Since $h$ is of class $C^2$,
\[h(T,Y)=h_T(0,Y)T+\frac12h_{TT}(0,Y)T^2(1+o(1)),\] and
$Th_T(T,Y)=h_T(0,Y)T+h_{TT}(0,Y)T^2(1+o(1))$. For $T=0$, we find
$w_1(0,Y)=\frac12 h_{TT}(0,Y)$. Since $h=0$ for $T=0$, we have
$h_{YY}(0,Y)=0$. Therefore, $w_1(0,Y)=\frac12 \Delta'h(0,Y)=-\frac12
\tilde{k}(0,Y)=-\frac12 k(0,Y)$. If $k=-2\Delta d$, we find
\[w_1(0,Y)=-\kappa(Y).
\]
\end{remark}
Let us now prove that $w_1\in C^{2+\alpha}_\sharp(\Omega'')$, and that
\begin{align*}
G : C^{\alpha}(\overline{\Omega''})&\To C^{2+\alpha}_\sharp(\Omega'')\\
k&\mapsto w_1,
\end{align*}
is a bounded operator.

First, let us transform the definition of $h$. Suppressing the $Y$
dependence, we have
\begin{gather*}
 h/T=\int_0^1h_T(T\sigma)\,d\sigma,
\\
 D(h/T)=\int_0^1T\sigma h_{TT}(T\sigma)\,d\sigma,
\end{gather*}
and finally,
\[ w_1=\int_0^1 \sigma h_{TT}(T\sigma)\,d\sigma,
\]
which proves that
\[
w_1\in C^\alpha(\overline{\Omega''}).
\] But we also have
\[ \Delta'(Dh)=(Th_T)_{TT}+(Th_T)_{YY}=D\Delta'
h+2h_{TT}=-D\tilde{k}+2h_{TT}\in C^\alpha(\overline{\Omega''}).
\]
Since $Dh=0$ for $T=0$ and $T=\theta$, and $Dh$ is bounded over
$\Omega''$, we find that
\[
Dh \text{ also is of class }C^{2+\alpha}(\overline{\Omega''}),
\]
using the usual Schauder estimates for this equation for $Dh$.
This proves
\[
T^2w_1=(D-1)h\in C^{2+\alpha}(\overline{\Omega''}).
\]

Since $Dw_1=T^{-2}(D-1)(D-2)h=T^{-2}[D(D-1)-2(D-1)]h=h_{TT}-2w_1$,
$(D+2)w_1$ is of class $C^\alpha$, and
\[
Tw_1\in C^{1+\alpha}(\overline{\Omega''}).
\]

Finally, let us show that $L_0w_1=k$.
\begin{align*}
L_0w_1 &=(D+2)(D-1)T^{-2}(D-1)h+(D-1)\pa^2_Y h\\
       &=T^{-2}D(D-3)(D-1)h+(D-1)\left\{-T^{-2}D(D-1)h-\tilde{k}\right\}\\
       &=T^{-2}D(D-1)(D-3)h-T^{-2}(D-3)D(D-1)h-(D-1)\tilde{k}\\
       &=k.
\end{align*}
This completes the proof.

\subsection{Solution of $Lw_0+2\Delta d=0$.}
\label{ssec:sol-w0}

We now treat equation $(L_0+L_1)w_0=-2\Delta d$ by a perturbation argument.
The previous section provides a bounded operator $G :
C^\alpha(\overline{\Omega''})\To C^{2+\alpha}_\sharp(\Omega'')$,
which is a right inverse for $L_0$.
We must now solve
\[ w_0=G[-2\Delta d]-G[L_1w_0].
\]
Since, by equation (\ref{eq:pert}), $w\mapsto G[L_1w]$ is a
contraction on $C^{2+\alpha}_\sharp(\Omega'')$ if $\theta$ is
sufficiently small, the result follows from the contraction
mapping principle.
\begin{remark}\label{rem:kappa2}
For $w_0\in C^{2+\alpha}_\sharp(\Omega'')$, it now follows from the
definition of $L_1$ that $L_1w_0$ vanishes for $T=Y=0$. It follows
from remark (\ref{rem:kappa}) that
$w_0(0,0)=\Delta d(0,0)=-\kappa(0)$. Theorem \ref{th:e5} shows that
$w(T,0)=w_0(T,0)+O(T\ln T)$, hence
\[ v(T,0)=2T-T^2(\kappa(0)+o(1)),
\]
which justifies the expansion (\ref{eq:exp}) in the introduction.
\end{remark}


\section{Construction of sub- and super-solutions
         and proof of theorem \ref{th:e5}}
\label{sec:e5}

We prove theorem \ref{th:e5}, using the information that $w$ and
$d\nabla w$ are bounded (theorems \ref{th:e2} and \ref{th:e3}).

Let $u_A=-\ln[2d+d^2w_A]$, where $w_A=w_0+A d\ln d$. This function
$u_A$ is, for $d<1$, an increasing function of $A$.
Taking $\Omega'$ smaller if necessary, we may assume that $w$ and
$w_A$ are bounded, and $|dw|$ and $|dw_A|$ are both less than one
over $\Omega'$; also, recall that $\pa\Omega'=\pa\Omega\cup\Gamma$.

Furthermore,
\[ L_0 (T\ln T)=(D+2)(D-1)T\ln T=T(D+3)D\ln T=3T,
\]
and $L_1(T\ln T)=O(T^2\ln T)$;
\[ Lw_A + (2+dw_A)\Delta d = AT(3+O(T\ln T))+2Tw_0\Delta d.
\]

Let us choose $A$ large enough and $\Omega'$ (i.e., the parameter
$\delta$) small enough so that
\[ u_{-A}\leq u \leq u_A
\]
on $\Gamma$, and
\begin{align*}
 Lw_{A} + (2+dw_A)\Delta d &\geq T\\
 Lw_{-A} + (2+dw_{-A})\Delta d &\leq -T.
\end{align*}
over $\Omega'$. We then find, by inspection of the expression for
$M_{w_A}w_A$, that
\[ (L-M_{w_A})w_{A} + 2\Delta d \geq T(1+\psi(A,T)),
\]
where $\psi(A,T)=O(T\ln T)$ for fixed $A$. We conclude from
(\ref{eq:ssol}) that $u_A$ is a super-solution of
(\ref{eq:liouville}) near the boundary if $A$ is large and
positive.

Similarly, $u_{-A}$ is a sub-solution near the boundary if $A$ is
large and negative.

Let us show that $w_A \leq w \leq w_{-A}$ over $\Omega'$.

\begin{lemma}
For any real $A$, $u-u_A=O(d)$ and $\nabla(u-u_A)=O(1)$ as $d\to
0$.
\end{lemma}
\begin{proof}
Since the function $t\mapsto \ln(2+t)$ has a bounded derivative
over $[-1,1]$.
\[
u-u_A=\ln(2+dw)-\ln(2+dw_A)=O(d(w-w_A))=O(d)
\]
since $w$ and $w_A$ are both bounded.

Next,
\begin{align*}
\nabla(u_A-u)&=\frac{\nabla(dw)}{2+dw}-\frac{\nabla(dw_A)}{2+dw_A}\\
             &=\left[\frac{w}{2+dw}-\frac{w}{2+dw_A}\right]\nabla d
               +\left[\frac{d\nabla w}{2+dw}-\frac{d\nabla
               w_A}{2+dw_A}\right]\\
             &=O(1)
\end{align*}
since $d\nabla w$ and $d\nabla w_A$ are both bounded. This
completes the proof.
\end{proof}
\begin{lemma}
Let $u_1$ and $u_2$ be respectively a sub- and a super-solution of
class $C^1(\Omega'\cup\Gamma)$ of equation (\ref{eq:liouville}) on
$\Omega'$. Assume that $u_1\leq u_2$ on $\Gamma$, and that
$(u_1-u_2)(x,y)=O(d)$ and $\nabla(u_1-u_2)(x,y)=O(1)$ as
$(x,y)\to\pa\Omega$. Then $u_1\leq u_2$ on $\Omega'$.
\end{lemma}
\begin{remark} This type of argument is taken from \cite{BBL}, see
also \cite{MV}.
\end{remark}
\begin{proof}
Let $\varphi$ be a smooth cut-off function equal to 1 if
$d>2\sigma$, zero if $d<\sigma$, and such that $0\leq \varphi\leq
1$ and $|\nabla\varphi|=O(1/\sigma)$. Testing the equation
$-\Delta(u_1-u_2)+4(e^{u_1}-e^{u_2})\leq 0$ with
$\varphi(u_1-u_2)_+$, which vanishes both on $\pa\Omega$ and on
$\Gamma$, and using the fact that $(u_1-u_2)(e^{u_1}-e^{u_2})\geq
0$, we find
\[
\int_{d>\sigma}\varphi|\nabla[(u_1-u_2)_+]|^2dx\,dy
+\int_{\sigma<d<2\sigma}(u_1-u_2)_+\nabla\varphi\cdot\nabla(u_1-u_2)dx\,dy
\leq 0.
\]
Consider now the second integral: it extends over the set where
$\sigma<d<2\sigma$, which has measure $O(\sigma)$; the integrand
on the other hand is, using the assumptions on $u_2-u_1$,
$O(\sigma)\times O(1/\sigma)=O(1)$. This second integral therefore
tends to zero with $\sigma$. It follows that
$\nabla[(u_1-u_2)_+]$, hence $(u_1-u_2)_+$, vanishes identically,
hence $u_1\leq u_2$, as desired.
\end{proof}
Applying this result to $u_{-A}$ and $u$, and then to $u$ and
$u_A$, we find that $u_{-A}\leq u \leq u_A$ on $\Omega'$.

This implies that $w_{A}\leq w \leq w_{-A}$, hence
\[ |w-w_0|\leq Ad\ln(1/d),  \]
QED.


%
\end{document}